\documentstyle[12pt]{article}
\newtheorem{theorem}{Theorem}[section]
\newtheorem{corollary}[theorem]{Corollary}

\newtheorem{example}[theorem]{Example}
\newtheorem{remark}[theorem]{Remark}
\newtheorem{lemma}[theorem]{Lemma}
\newtheorem{proposition}[theorem]{Proposition}
\title{More on sg-compact spaces\thanks{1991 Math.\ Subject
Classification --- Primary: 54D30, 54A05; Secondary: 54H05,
54G99. \newline Key words and phrases --- sg-compact,
semi-compact, $C_2$-space, semi-open set, sg-open set, hsg-closed
sets. \newline Research supported partially by the Ella and Georg
Ehrnrooth Foundation at Merita Bank, Finland.}}
\author{Julian Dontchev\\Department of Mathematics\\University
of Helsinki\\PL 4, Yliopistonkatu 15\\00014 Helsinki 10\\Finland
\and Maximilian Ganster\\Department of Mathematics\\Graz
University of Technology\\Steyrergasse 30\\A-8010 Graz\\Austria}
\date{}
\begin{document}
\baselineskip=20pt plus 1pt minus 1pt
\maketitle
\begin{abstract}
The aim of this paper is to continue the study of sg-compact
spaces, a topological notion much stronger than hereditary
compactness. We investigate the relations between sg-compact and
$C_2$-spaces and the interrelations to hereditarily sg-closed
sets.
\end{abstract}

\section{Introduction}\label{s1}

In 1995, sg-compact spaces were introduced independently by
Caldas \cite{CC1} and by Devi, Balachandran and Maki \cite{DBM1}.
A topological space $(X,\tau)$ is called {\em sg-compact}
\cite{CC1} if every cover of $X$ by sg-open sets has a finite
subcover. In \cite{DBM1}, the term {\em $SGO$-compact} is used.

Recall that a subset $A$ of a topological space $(X,\tau)$ is
called {\em sg-open} \cite{BL1} if every semi-closed subset of
$A$ is included in the semi-interior of $A$. A set $A$ is called
{\em semi-open} if $A \subseteq \overline{{\rm Int} A}$ and {\em
semi-closed} if ${\rm Int} \overline{A} \subseteq A$. The {\em
semi-interior} of $A$, denoted by ${\rm sInt}(A)$, is the union
of
all semi-open subsets of $A$ while the {\em semi-closure} of $A$,
denoted by ${\rm sCl}(A)$, is the intersection of all semi-closed
supersets of $A$. It is well known that ${\rm sInt}(A)$ =
$A \cap \overline{{\rm Int}A}$ and ${\rm sCl} (A)$ = $A \cup 
{\rm Int}\overline {A}$ .

Every topological space $(X,\tau)$ has a unique decomposition
into two sets $X_1$ and $X_2$, where $X_1 = \{ x \in X \colon \{
x \}$ is nowhere dense$\}$ and $X_2 = \{ x \in X \colon \{ x \}$
is locally dense$\}$. This decomposition follows from a result
of
Jankovi\'{c} and Reilly \cite[Lemma 2]{JR1}. Recall that a set
$A$
is said to be {\em locally dense} \cite{CM1} (= {\em preopen})
if 
$A \subseteq {\rm Int} \overline{A}$.\vspace {6 mm}

It is a fact that a subset $A$ of $X$ is sg-closed (= its
complement is sg-open) if and only if $X_1 \cap {\rm sCl}(A)
\subseteq A$ \cite{DM1}, or equivalently  if and only if $X_1 
\cap {\rm Int} \overline{A} \subseteq A$. By taking complements 
one easily observes that $A$ is sg-open if and only if $A \cap 
X_1 \subseteq {\rm sInt}(A)$. Hence every subset of $X_2$ is
sg-open.

\section{Sg-compact spaces}\label{s2}

Let $A$ be a sg-closed subset of a topological space $(X,\tau)$.
If every subset of $A$ is also sg-closed in $(X,\tau)$, then $A$
will be called {\em hereditarily sg-closed} (= hsg-closed).
Observe that every nowhere dense subset is hsg-closed but not
vice versa.

\begin{proposition}\label{p1}
For a subset $A$ of a topological space $(X,\tau)$ the following
conditions are equivalent:

{\rm (1)} $A$ is hsg-closed.

{\rm (2)} $X_1 \cap {\rm Int} \overline{A} = \emptyset$.
\end{proposition}

{\em Proof.} (1) $\Rightarrow$ (2) Suppose that there exits $x
\in X_1 \cap {\rm Int} \overline{A}$. Let $V_x$ be an open set
such that $V_x \subseteq \overline{A}$ and let $B = A \setminus
\{ x \}$. Since $B$ is sg-closed, i.e.\ $X_1 \cap {\rm sCl}(B)
\subseteq B$, we have $x \not\in {\rm sCl}(B)$, hence $x \not\in
{\rm Int} \overline{B}$, and thus $x \in \overline{X \setminus
\overline{B}}$. If $H = V_x \cap (X \setminus \overline{B})$,
then $H$ is nonempty and open with $H \subseteq \overline{A}$ and
$H \cap B = \emptyset$ and so $H \cap A = \{ x \}$. Hence
$\emptyset \not= H = H \cap \overline{A} \subseteq \overline{H
\cap A} \subseteq \overline{\{x\}}$, i.e.\ ${\rm Int}
\overline{\{x\}} \not= \emptyset$. Thus $x \in X_2$, a
contradiction.

(2) $\Rightarrow$ (1) Let $B \subseteq A$. Then ${\rm Int}
\overline{B} \subseteq {\rm Int} \overline{A}$ and $X_1 \cap {\rm
Int} \overline{B} = \emptyset$, i.e.\ $B$ is sg-closed. $\Box$
\vspace {6 mm}

We will call a topological space $(X,\tau)$  a {\em $C_2$-space}
\cite{G1} (resp.\ {\em $C_3$-space}) if every nowhere dense
(resp.\ hsg-closed) set is finite. Clearly every $C_3$-space is
a $C_2$-space. Also, a topological space $(X,\tau)$ is indiscrete
if and only if every subset of $X$ is hsg-closed (since in that
case $X_1 =\emptyset$).

Following Hodel \cite{H1}, we say that a {\em cellular family}
in a topological space $(X,\tau)$ is a collection of nonempty,
pairwise disjoint open sets. The following result reveals an
interesting property of $C_2$-spaces.

\begin{lemma}\label{l1}
Let $(X,\tau)$ be a $C_2$-space. Then every infinite cellular
family has an infinite subfamily whose union is contained in
$X_2$.
\end{lemma}

{\em Proof.} Let $\{ U_i \colon i \in {\bf N} \}$ be a cellular
family. Suppose that for infinitely many $i \in {\bf N}$ we have
$U_i \cap X_1 \not= \emptyset$. Without loss of generality
we may assume that $U_i \cap X_1 \not= \emptyset$ for each $i \in
{\bf N}$. Now pick $x_i \in U_i \cap X_1$ for each $i \in {\bf
N}$ and
partition $\bf N$ into infinitely many
disjoint infinite sets, ${\bf N} = \cup_{k \in {\bf N}} {\bf
N}_{k}$. Let $A_k = \{ x_i \colon i \in {\bf N}_{k} \}$. Since
$A_k \cap ( \cup_{i \not\in {\bf N}_{k}} U_{i}) = \emptyset$ and
$A_k \subseteq \cup_{i \in {\bf N}_{k}} U_{i}$ for each $k$,
it is easily checked that $\{ {\rm Int} \overline{A_{k}} \colon
k \in {\bf N} \}$ is a disjoint family of open sets. Since $X$
is a $C_2$-space, $A_k$ cannot be nowhere dense and so, for
each $k$, there exists $p_k \in {\rm Int} \overline{A_{k}}$ and
the $p_k$'s are pairwise distinct. Also, since $X$ is $C_2$,
$\overline{\cup_{i \in {\bf N}} U_{i}} = \cup_{i \in {\bf N}}
(U_{i}) \cup F$, where $F$ is finite. Since $p_k \in
\overline{\cup_{i \in {\bf N}} U_{i}}$ for each $k$, there exists
$k_0$ such that $p_k \in \cup_{i \in {\bf N}} U_{i}$ for $k \geq
k_0$, and
since ${\rm Int} \overline{A_{k}} \cap (\cup_{i \not\in {\bf
N}_{k}} U_{i})
= \emptyset$, we have $p_k \in \cup_{i \in {\bf N}_{k}} U_{i}$
for
$k \geq k_0$. Now, for each $k \geq k_0$ pick $i_k \in {\bf
N}_{k}$ such that $p_k \in U_{i_{k}}$, and so $p_k \in W =
U_{i_{k}} \cap
{\rm Int} \overline{A_{k}}$. Thus $\emptyset \not= W \subseteq
U_{i_{k}} \cap \overline{A_{k}} \subseteq \overline{U_{i_{k}}
\cap A_{k}} = \overline{\{ x_{i_{k}} \}}$. Hence $\{ x_{i_{k}}
\}$ is locally dense, a contradiction. This shows that only for
finitely many $i \in {\bf N}$ we have $U_i \cap X_1 \not=
\emptyset$. Thus the claim is proved. $\Box$ \vspace{6mm}

The {\em $\alpha$-topology} \cite{Nj1} on a topological space
$(X,\tau)$ is the collection of all sets of the form $U \setminus
N$, where $U \in \tau$ and $N$ is nowhere dense in $(X,\tau)$.
Recall that topological spaces whose $\alpha$-topologies are
hereditarily compact have been shown to be {\em semi-compact}
\cite{GJR1}.
The original definition of semi-compactness is in terms of
semi-open sets and
is due to Dorsett \cite{D1}. By definition a topological space
$(X,\tau)$ is called {\em semi-compact} \cite{D1} if every cover
of $X$ by semi-open sets has a finite subcover.

\begin{remark}\label{r1}
{\em (i) The 1-point-compactification of an infinite discrete
space is a $C_2$-space having an infinite cellular family.

(ii) \cite{G1} A topological space $(X,\tau)$ is semi-compact
if and only if $X$ is a $C_2$-space and every cellular family is
finite.

(iii) \cite{HD1} Every subspace of a semi-compact space is
semi-compact (as a subspace).}
\end{remark}

\begin{lemma}\label{l2}
{\rm (i)} Every $C_3$-space $(X,\tau)$ is semi-compact.

{\rm (ii)} Every sg-compact space is semi-compact.
\end{lemma}

{\em Proof.} (i) All $C_3$-spaces are $C_2$-spaces. Thus in the
notion of Remark~\ref{r1} (ii) above we need to show that every
cellular family in $X$ is finite. Suppose that there exists an
infinite cellular family $\{ U_i \colon i \in {\bf N} \}$. For
each $i \in {\bf N}$ pick $x_i \in U_i$ and, as before, partition
${\bf N} = \cup_{k} {\bf N}_{k}$ and set $A_k = \{ x_i \colon i
\in {\bf N}_{k} \}$. Since $X$ is a $C_2$-space, $\{ {\rm
Int} \overline{A_{k}} \colon k \in {\bf N} \}$ is a cellular
family. By Lemma~\ref{l1}, there is a $k \in {\bf N}$ such that
${\rm Int} \overline{A_{k}} \subseteq X_2$. Since $A_k$ is not
hsg-closed, we must have $X_1 \cap {\rm Int} \overline{A_{k}}
\not=
\emptyset$, a contradiction. So, every cellular family in $X$ is
finite and consequently $(X,\tau)$ is semi-compact.

(ii) is obvious since every semi-open set is sg-open. $\Box$

\begin{remark}
{\em (i) It is known that sg-open sets are $\beta$-open, i.e.\
they are dense in some regular closed subspace \cite{JD1}. Note
that $\beta$-compact spaces, i.e.\ the spaces in which every
cover by $\beta$-open sets has a finite subcover are finite
\cite{G2}. However, one can easily find an example of an infinite
sg-compact space -- the real line with the cofinite topology is
such a space.

(ii) In semi-$T_D$-spaces the concepts of sg-compactness and
semi-compactness coincide. Recall that a topological space
$(X,\tau)$ is called a {\em semi-$T_D$-space} \cite{JR1} if each
singleton is either open or nowhere dense, i.e.\ if every
sg-closed set is semi-closed.}
\end{remark}

\begin{theorem}\label{t1}
For a topological space $(X,\tau)$ the following conditions are
equivalent:

{\rm (1)} $X$ is sg-compact.

{\rm (2)} $X$ is a $C_3$-space.
\end{theorem}

{\em Proof.} (1) $\Rightarrow$ (2) Suppose that there exists an
infinite hsg-closed set $A$ and set $B = X \setminus A$. Observe
that for each $x \in A$, the set $B \cup \{ x \}$ is sg-open in
$X$. Thus $\{ B \cup \{ x \} \colon x \in A \}$ is a sg-open
cover of $X$ with no finite subcover. Thus $(X,\tau)$ is $C_3$.

(2) $\Rightarrow$ (1) Let $X = \cup_{i \in I} A_i$, where each
$A_i$ is sg-open. Let $S_i = {\rm sInt}(A_{i})$ for each $i \in
I$ and let $S = \cup_{i \in I} S_i$. Then $S$ is a semi-open
subset of $X$ and each $S_i$ is a semi-open subset of
$(S,\tau|S)$. Since $X$ is a $C_3$-space,  $(X,\tau)$ is
semi-compact and hence $(S,\tau|S)$ is a semi-compact subspace
of $X$ (by Remark~\ref{r1} (iii)). So we may say that
 $S = S_{i_{1}} \cup \ldots \cup S_{i_{k}}$.
Since $A_i$ is sg-open, we have $X_1 \cap
A_i \subseteq S_i$ for each index $i$ and so $X_1 = X_1 \cap
(\cup A_{i}) \subseteq X_1 \cap S \subseteq S_{i_{1}} \cup \ldots
\cup S_{i_{k}} = S$. Hence $X \setminus S$ is semi-closed and $X
\setminus S \subseteq X_2$. Since ${\rm Int} \overline{(X
\setminus S)} \subseteq X \setminus S \subseteq X_2$, we conclude
that $X \setminus S$ is hsg-closed and thus finite. This shows
that $X
= S_{i_{1}} \cup \ldots \cup S_{i_{k}} \cup (X \setminus S) =
A_{i_{1}} \cup \ldots \cup A_{i_{k}} \cup F$, where $F$ is
finite, i.e.\ $(X,\tau)$ is sg-compact. $\Box$

\begin{remark}\label{r2}
{\em (i) If $X_1 = X$, then $(X,\tau)$ is sg-compact if and only
if $(X,\tau)$ is semi-compact. Observe that in this case
sg-closedness and semi-closedness coincide.

(ii) Every infinite set endowed with the cofinite topology is
(hereditarily) sg-compact.}
\end{remark}

It is known that an arbitrary intersection of sg-closed sets is
also
an sg-closed set \cite{DM1}. The following result provides an
answer to the question about the additivity of sg-closed sets.

\begin{proposition}\label{r3}
{\rm (i)} If $A$ is sg-closed and $B$ is closed, then $A \cup B$
is also sg-closed.

{\rm (ii)} The intersection of a sg-open and an open set is
always sg-open.

{\rm (iii)} The union of a sg-closed and a semi-closed set need
not be sg-closed, in particular, even finite union of sg-closed
sets need not be sg-closed.
\end{proposition}

{\em Proof.} (i) Let $A \cup B \subseteq U$, where $U$ is
semi-open. Since $A$ is sg-closed, we have ${\rm sCl}(A \cup B)
=
(A \cup B) \cup {\rm Int} (\overline{A \cup B}) \subseteq U \cup
{\rm Int} (\overline{A} \cup B) \subseteq U \cup ({\rm Int}
\overline{A}
\cup B) \subseteq U \cup (U \cup B) = U$.

(ii) follows from (i).

(iii) Let $X = \{ a,b,c,d \}$, $\tau = \{ \emptyset, \{ a \}, \{
b \}, \{ a,b \}, X \}$. Note that the two sets $A = \{ a \}$ and
$B = \{ b \}$ are semi-closed but their union $\{ a,b \}$ is not
sg-closed. $\Box$ \vspace{6mm}

Theorem 3 from \cite{BL1} states that if $B \subseteq A \subseteq
(X,\tau)$ and $A$ is open and sg-closed, then $B$ is sg-closed
in the subspace $A$ if and only if $B$ is sg-closed in $X$. Since
a subset is regular open if and only if it is $\alpha$-open and
sg-closed \cite{DP1}, by using Proposition~\ref{r3},
we obtain the following result:

\begin{proposition}\label{p4}
Let $R$ be a regular open subset of a topological space
$(X,\tau)$. If $A \subseteq R$ and $A$ is sg-open in
$(R,\tau|R)$, then $A$ is sg-open in $X$. $\Box$
\end{proposition}

{\em Proof.} Since $B = R \setminus A$ is sg-closed in
$(R,\tau|R)$, $B$ is sg-closed in $X$ by \cite[Theorem 3]{BL1}. 
Thus $X \setminus B$ is sg-open in $X$ and by
Proposition~\ref{r3}
(ii), $R \cap (X \setminus B) = A$ is sg-open in $X$. $\Box$
\vspace{6mm}

Recall that a subset $A$ of a topological space $(X,\tau)$ is
called {\em $\delta$-open} \cite{V1} if $A$ is a union of regular
open sets. The collection of all $\delta$-open subsets of a
topological space $(X,\tau)$ forms the so called {\em
semi-regularization topology}.

\begin{corollary}\label{c2}
If $A \subseteq B \subseteq (X,\tau)$ such that $B$ is
$\delta$-open in $X$ and $A$ is sg-open in $B$, then $A$ is
sg-open in $X$.
\end{corollary}

{\em Proof.} Let $B = \cup_{i \in I} B_i$, where each $B_i$ is
regular open in $(X,\tau)$. Clearly, each $B_i$ is regular open
also in $(B,\tau|B)$. By Proposition~\ref{r3} (ii), $A \cap B_i$
is sg-open in $(B,\tau|B)$ for each $i \in I$. In the notion of
Proposition~\ref{p4}, $B \setminus (A \cap B_{i})$ is sg-closed
in $(X,\tau)$ for each $i \in I$. Hence $X \setminus (B \setminus
(A \cap B_{i})) = (A \cap B_{i}) \cup (X \setminus B)$ is sg-open
in $(X,\tau)$. Again by Proposition~\ref{r3} (ii), $B \cap ((A
\cap B_i) \cup (X \setminus B)) = A \cap B_i$ is sg-open in
$(X,\tau)$. Since any union of sg-open sets is always sg-open,
we have
$A = \cup_{i \in I} (A \cap B_i)$ is sg-open in $(X,\tau)$.
$\Box$

\begin{proposition}\label{p2}
Every $\delta$-open subset of a sg-compact space $(X,\tau)$ is
sg-compact, in particular, sg-compactness is hereditary with
respect to regular open sets.
\end{proposition}

{\em Proof.} Let $A \subseteq X$ be $\delta$-open. If $\{ U_i
\colon i \in I \}$ is a sg-open cover of $(S,\tau|S)$, then by
Corollary~\ref{c2}, each $U_i$ is sg-open in $X$. Then, $\{ U_i
\colon i \in I \}$ along with $X \setminus A$ forms a sg-open
cover of $X$. Since $X$ is sg-compact, there exists a finite $F
\subseteq I$ such that $\{ U_i \colon i \in F \}$ covers $A$.
$\Box$

\begin{example}\label{e1}
{\em Let $A$ be an infinite set with $p \not\in A$. Let $X = A
\cup \{ p \}$ and $\tau = \{ \emptyset, A, X \}$.

(i) Clearly, $X_1 = \{ p \}$, $X_2 = A$ and for each infinite $B
\subseteq X$, we have $\overline{B} = X$. Hence $X_1 \cap {\rm
Int} \overline{B} \not= \emptyset$, so $B$ is not hsg-closed.
Thus $(X,\tau)$ is a $C_3$-space, so sg-compact. But the open
subspace $A$ is an infinite indiscrete space which is not
sg-compact. This shows that (1) hereditary sg-compactness is
a strictly stronger concept than sg-compactness and (2) in
Proposition~\ref{p2} '$\delta$-open' cannot be replaced with
'open'.

(ii) Observe that $X \times X$ contains an infinite nowhere dense
subset, namely  $X \times X \setminus A \times A$. This shows
that even
the finite product of two sg-compact spaces need not be
sg-compact, not even a $C_2$-space.

(iii) \cite{MBD1} If the nonempty product of two spaces is
sg-compact $T_{gs}$-space (see \cite{MBD1}), then each
factor space is sg-compact.}
\end{example}

Recall that a function $f \colon (X,\tau) \rightarrow
(Y,\sigma)$ is called {\em pre-sg-continuous} \cite{N1} if
$f^{-1}(F)$ is sg-closed in $X$ for every semi-closed subset $F
\subseteq Y$.

\begin{proposition}
{\rm (i)} The property 'sg-compact' is topological.

{\rm (ii)} Pre-sg-continuous images of sg-compact spaces are
semi-compact. $\Box$
\end{proposition}

\
\
E-mail: {\tt dontchev@cc.helsinki.fi}, {\tt
ganster@weyl.math.tu-graz.ac.at}
\
\

\baselineskip=12pt

\begin{thebibliography}{18}\frenchspacing


\bibitem{BL1} {P. Bhattacharyya and B.K. Lahiri},
{Semi-generalized closed sets in topology}, {\em Indian J.
Math.}, {\bf 29} (3) (1987), 375--382.

\bibitem{CC1} {M.C. Caldas}, {Semi-generalized continuous maps
in topological spaces}, {\em Portug. Math.}, {\bf 52} (4) (1995),
399--407.

\bibitem{CM1} {H.H. Corson and E. Michael}, {Metrizability of
certain countable unions}, {\em Illinois J. Math.}, {\bf 8}
(1964), 351--360.

\bibitem{DBM1} {R. Devi, K. Balachandran and H. Maki},
{Semi-generalized homeomorphisms and generalized
semi-homeomorphisms in topological spaces}, {\em Indian J. Pure
Appl. Math.}, {\bf 26} (3) (1995), 271--284.

\bibitem{JD1} {J. Dontchev}, {On some separation axioms
associated with the $\alpha$-topology}, {\em Mem. Fac. Sci. Kochi
Univ. Ser. A Math.}, {\bf 18} (1997), 31--35.

\bibitem{DM1} {J. Dontchev and H. Maki}, {On sg-closed sets and
semi-$\lambda$-closed sets}, {\em Questions Answers Gen.
Topology}, {\bf 15} (2) (1997), to appear.

\bibitem{DP1} {J. Dontchev and M. Przemski}, {On the various
decompositions of continuous and some weakly continuous
functions}, {\em Acta Math. Hungar.}, {\bf 71} (1-2) (1996),
109--120.

\bibitem{D1} {Ch. Dorsett}, {Semi-compact $R_1$ and product
spaces}, {\em Bull. Malaysian Math. Soc.}, {\bf 3} (2) (1980),
15--19.

\bibitem{G1} {M. Ganster}, {Some remarks on strongly compact
spaces and semi-compact spaces}, {\em Bull. Malaysia Math. Soc.},
{\bf 10} (2) (1987), 67--81.

\bibitem{G2} {M. Ganster}, {Every $\beta$-compact space is
finite}, {\em Bull. Calcutta Math. Soc.}, {\bf 84} (1992),
287--288.

\bibitem{GJR1} {M. Ganster, D.S. Jankovi\'{c} and I.L. Reilly},
{On compactness with respect to semi-open sets}, {\em Comment.
Math. Univ. Carolinae}, {\bf 31} (1) (1990), 37--39.

\bibitem{HD1} {F. Hanna and Ch. Dorsett}, {Semicompactness}, {\em
Questions Answers Gen. Topology}, {\bf 2} (1) (1984), 38--47.

\bibitem{JR1} {D. Jankovi\'{c} and I. Reilly}, {On semiseparation
properties}, {\em Indian J. Pure Appl. Math.}, {\bf 16} (9)
(1985), 957--964.

\bibitem{H1} {R. Hodel}, {Cardinal Functions I}, Handbook of
Set-Theoretic Topology, North Holland (1987).

\bibitem{MBD1} {H. Maki, K. Balachandran and R. Devi}, {Remarks
on semi-generalized closed sets and generalized semi-closed
sets}, {\em Kyungpook Math. J.}, {\bf 36} (1996), 155-163.

\bibitem{Nj1} {O. Nj{\aa}stad}, {On some classes of nearly open
sets}, {\em Pacific J. Math.}, {\bf 15} (1965), 961--970.

\bibitem{N1} {T. Noiri}, {Semi-normal spaces and some functions},
{\em Acta Math. Hungar.}, {\bf 65} (3) (1994), 305--311.

\bibitem{V1} {N.V. Veli\v{c}ko}, {$H$-closed topological spaces},
{\em Amer. Math. Soc. Transl.}, {78} (1968), 103--118.


\end{thebibliography}
\end{document}